\documentclass[preprint,12pt]{elsarticle}

\usepackage{natbib}
\usepackage{color}
\usepackage{graphicx}
\usepackage{amssymb}
\newtheorem{theorem}{Theorem}[section]
\newtheorem{definition}{Definition}[section]

\newtheorem{lemma}{Lemma}[section]

\journal{arXiv.org}

\begin{document}

\begin{frontmatter}

\title{Cyber-Physical Systems as General Distributed Parameter Systems:  Three Types of Fractional Order Models  and Emerging Research Opportunities}

\author[dhu,ucmerced]{Fudong Ge}
\ead{gefd2011@gmail.com}
\author[ucmerced]{YangQuan Chen\corref{cor}}
\ead{ychen53@ucmerced.edu}
\author[dhutwo]{Chunhai Kou}
\ead{kouchunhai@dhu.edu.cn}
\cortext[cor]{Corresponding author}

\address[dhu]{College of Information Science and Technology, Donghua University, Shanghai 201620, PR China}
\address[ucmerced]{Mechatronics, Embedded Systems and Automation Lab, University of California, Merced, CA 95343, USA}
\address[dhutwo]{Department of Applied Mathematics, Donghua University, Shanghai 201620, PR China}

\begin{abstract}
Cyber-physical systems (CPSs) are man-made complex systems coupled with natural processes that, as a whole,  should be described by distributed
parameter systems (DPSs) in general forms. This paper presents three such general models for generalized DPSs that can
be used to characterize complex CPSs. These three different types of fractional operators based DPS models  are: fractional Laplacian operator,
fractional power of operator or fractional derivative. This research investigation is motivated by many fractional order models describing natural,
physical, and anomalous phenomena, such as sub-diffusion process or super-diffusion process.  The  relationships  among  these three  different
operators are explored and explained. Several  potential future research opportunities are then articulated followed by  some conclusions and remarks.
\end{abstract}
\par

\end{frontmatter}

\section{Introduction}
\setcounter{section}{1}\setcounter{equation}{0}
It is well known that the CPSs with integrated computational and physical processes can be regarded as
a new generation of control systems and can interact with humans through many new modalities \cite{CPSlee2008}. The objective of CPS
is to develop new science and engineering methods in which cyber and physical designs are compatible, synergistic, and integrated
at all scales. Besides, as we all know, the DPSs  can be  used to well characterize those cyber-physical
process \cite{CPSsong2009,CPStricaud2011} and the actions and measurements of the system studied are better described by utilizing
the actuators and sensors, which was first introduced by El Jai and Pritchard in \cite{AEIJAI} and mainly focused on the
locations, number and spatial distributions of the actuators and sensors.

Moreover, in the past several decades, fractional calculus has shown great potential in science and engineering applications and some phenomena such as
self-similarity, non-stationary, non-Gaussian process and short or long memory process are all closely related to fractional calculus \cite{Podlubny,Kilbas,Klimek}.
It is now widely believed that, using fractional calculus in modeling can better capture the complex dynamics of natural and man-made systems, and fractional
order controls can offer better performance not achievable before using integer order controls \cite{Torvik,Mandelbrot}.

Motivated by the above arguments, in this paper, let $\Omega$ be an open bounded subset of $\mathbf{R}^n$ with smooth boundary
$\partial\Omega$ and we  consider the following fractional DPSs:
\begin{equation}\label{problem}
z_t(x,t)+Az(x,t)=u(x,t)~\mbox{ in } \Omega \times [0,b],
\end{equation}
where $b>0$ is a given constant, $u$ is the control input depending on the number and the structure of actuators and
$A$ may be a fractional Laplacian operator, a fractional power of operator or a fractional derivative.

The contribution of this present paper is to analyze the relationship among the fractional Laplacian operator, fractional power of
operator and fractional derivative and try to explore the opportunities and research challenges related to the fractional order DPSs
emerging at the same time. To the best of our knowledge, no result is available on this topic.  We hope that the results here could
provide some insights into the control theory of this field and be used in real-life applications.

The rest of the paper is organized as follows. The relationship among fractional Laplacian operator, fractional power of operator
and  fractional derivative are explored in Sec.~\ref{sec2}. In Sec.~\ref{sec3}, the emerging research opportunities of the fractional
order DPSs with those three operators are discussed. Several conclusions and remarks of this paper are given in the last section.

\section{Three different types of operators}
\label{sec2}

\setcounter{section}{2}\setcounter{equation}{0}
In this section, we shall introduce some basic relationships among the fractional Laplacian operator, fractional power of operator and fractional derivative. For further information, we refer the readers to papers from \cite{lv2015} to [36]
%
in the reference section of the present paper 
and the references cited  therein.

\subsection{Fractional Laplacian operator and fractional power of operator}

This subsection is devoted to the difference between the fractional Laplacian operator and fractional power of operator. For more details,
please see \cite{lv2015,lv2014,lv2013,di2012,du2013,applebaum2009} and their cited references.

Let us denote   $(-\triangle)^{\alpha/2}$ the nonlocal operator (also called the fractional Laplacian operator) defined
pointwise by the following Cauchy principal value integral
\begin{eqnarray}\label{cauchyd}
(-\triangle)^{\alpha/2} f(x)=C_{\alpha} P.V. \int_{\mathbf{R}}{\frac{f(x)-f(y)}{|x-y|^{1+\alpha}}}dy,~~0<\alpha<2,
\end{eqnarray}
where $C_{\alpha}=\frac{2^\alpha\Gamma(1/2+\alpha/2)}{\sqrt{\pi} \left|\Gamma\left(-\alpha/2\right)\right|}$ is a constant dependent on the order $\alpha.$
Obviously, we see that the fractional Laplacian $\left(-\triangle\right)^{\alpha/2}$ is a nonlocal operator which depends on the
parameter $\alpha$ and recovers the usual Laplacian as $\alpha\to 2$.
For more information about the fractional Laplacian operator, see \cite{du2012,choi2014,barrios2012,micu2006,fattorini1971}
and the references cited  therein.
Now we have the following result.
\begin{theorem}\cite{Kwasnicki2012}
Suppose that $(-\triangle)^{\alpha/2}$ is defined in $L^2(0, l)$ for $\alpha\in (0, 2)$. Then,  the eigenvalues of
the following  spectral problem
\begin{eqnarray}
(-\triangle)^{\alpha/2} \xi(x)=\lambda \xi (x),~~x\in (0,l),
\end{eqnarray}
where $\xi\in L^2(0, l)$ is extended to all $\mathbf{R}$ by $0$ is
\begin{eqnarray}
\lambda_n=\left(\frac{n\pi}{l}-\frac{(2-\alpha)\pi}{4l}\right)^\alpha+O\left(\frac{1}{n}\right)
\end{eqnarray}
satisfying
\[0<\lambda_1<\lambda_2\leq \cdots\leq \lambda_i\leq \cdots.\]
Moreover,  the corresponding eigenfunctions $\xi_n$ of $\lambda_n,$ after normalization, form a complete orthonormal basis in $L^2(0,l).$
\end{theorem}

Note that the constant in the error term  $O\left(\frac{1}{n}\right)$ tends to zero when $\alpha$ approaches 2 and in the particular case, when $\alpha=2$, we see that
$\lambda_n=\left(n\pi/l\right)^2$ without the error term.

However,  for a positive operator $A$ on bounded domain $[0,l]$. Suppose $0<\lambda_1\leq \lambda_2\leq  \cdots\leq \lambda_n\leq \cdots$
are the eigenvalues of $A$,  $\{\xi_1, \xi_2,\cdots,\xi_n, \cdots\}$
are the corresponding eigenfunctions and $\xi_n~(i=1,2,\cdots)$ form an orthonormal basis of $L^2(0,l)$. Let $(\cdot,\cdot)$  be the inner product of $L^2(0,l).$
We define the fractional power of operator $A$ as follows:
\begin{eqnarray}
A^\beta f(x)=\sum\limits_{n=1}^\infty{\lambda_n^\beta(\xi_n,f)\xi_n(x)},~~f\in L^2(0,l),
\end{eqnarray}
Then $\lambda_n^\beta~(i=1,2,\cdots)$ are the eigenvalues of $A^\beta$.
This implies that the two operators are different.

Moreover, the work spaces of the two operators (fractional Laplacian operator and fractional power of operator $A$) are different. Before   stating our main results, we first introduce two Banach spaces,
which are specified in \cite{lv2014,di2012,du2013,du2012}.

For $\Omega\subseteq \mathbf{R}^n$ is a bounded domain, $s\in (0,1)$ and $p\in [1,\infty)$, we define the classical Sobolev space $W^{s,p}(\Omega)$ as follows \cite{di2012}:
\begin{eqnarray}
W^{s,p}(\Omega):=\left\{  f\in L^p(\Omega): \frac{f(x)-f(y)}{|x-y|^{\frac{n}{p}+s}}\in L^p(\Omega \times \Omega)\right\}
\end{eqnarray}
endowed with the natural norm
\begin{eqnarray}
\|f\|_{W^{s,p}(\Omega)}:=\left( \int_\Omega{|f(x)|^p}dx+\int_\Omega\int_\Omega{\frac{|f(x)-f(y)|^p}{|x-y|^{n+sp}}}dxdy\right)^{1/p}
\end{eqnarray}
is an intermediary Banach space between $L^p(\Omega)$ and  $W^{1,p}(\Omega)$.
When a non-integer $s>1$, let $s=m+\sigma$ with  $m\in \mathbf{N}$ and  $\sigma\in (0,1)$.
In this case, let $D^\beta f$ with $|\beta|=m$  be the distributed derivative of $f$, then the classical Sobolev  space $W^{s,p}(\Omega)$  defined by
\begin{eqnarray}
W^{s,p}(\Omega):=\left\{  f\in W^{m,p}(\Omega): D^\beta f\in W^{\sigma,p}(\Omega) \mbox{ for all } \beta \mbox{ such that  }|\beta|=m\right\}
\end{eqnarray}
with respect to the norm
\begin{eqnarray}
\|f\|_{W^{s,p}(\Omega)}:=\left(\|f\|^p_{W^{m,p}(\Omega)}+\|D^\beta f\|^p_{W^{\sigma,p}(\Omega)}\right)^{1/p}
\end{eqnarray}
is a Banach space.  Clearly, if $s = m$ is an integer, the space $W^{s,p}(\Omega)$ coincides with the Sobolev space $W^{m,p}(\Omega)$.

Besides, let $\rho(x)\sim\frac{1}{\delta^\alpha(x)}$ with $\delta(x)=\mbox{dist}(x,\Omega^c)$. Define another space as follows:
\begin{eqnarray}
W^{s,p}_\rho(\Omega):=\left\{  f\in W^{s,p}(\Omega) : \rho(x)f(x)\in L^p(\Omega)\right\}
\end{eqnarray}
with the norm
\begin{eqnarray}
\|f\|_{W^{s,p}_\rho(\Omega)}:=\left( \int_\Omega{|\rho(x)f(x)|^p}dx+\int_\Omega\int_\Omega{\frac{|f(x)-f(y)|^p}{|x-y|^{n+sp}}}dxdy\right)^{1/p}.
\end{eqnarray}
Actually, $W^{s,p}_\rho(\Omega)$ is called nonlocal Sobolev space and we have $W^{s,p}_\rho(\Omega) \subseteq W^{s,p}(\Omega)$ \cite{du2013,du2012}.

By the Remark 2.1 in \cite{lv2014},
for the fractional power of operator, we take the classical fractional Sobolev space as its work space. But for fractional Laplacian operator,
we must take the nonlocal Sobolev space as its work space, which can be regarded as the weighted fractional Sobolev space. More precisely,
for any element  $f\in W^{s,p}(\Omega)$, since $-\triangle$ is a local operator and  we do not know how $f(x)$ approaches $0$ when $x \to \partial \Omega$.
Even if considering the space $W_0^{s,p}(\Omega);=\{f\in W^{s,p}(\Omega): f|_{\partial\Omega}=0\}$, we only know that the function
$f= 0$ on the boundary and we do not know how $f$ approaches 0.  However, for the fractional Laplacian operator, it is  a nonlocal operator and it
is defined in the whole space. So, it provides information about how $f$ approaches $0$ as $x \to \partial \Omega$. In fact, from the definition of nonlocal Sobolev space,
we know that $\frac{f(x)}{\delta^\alpha(x)}\to 0$  as $x \to \partial \Omega$, which dictates how $f$ approaches $0$ near boundary. It coincides with the result of the Theorem
1.2 in \cite{ros2014}. Thus, this is a significant difference between the fractional power of operator $A=-\triangle$ and the fractional Laplacian operator.


Besides, it is well known that $\|\triangle f\|_{W^{s,p}(\Omega)}=\| f\|_{W^{s+2,p}(\Omega)}$. But for the fractional Laplacian operator $(-\triangle)^\alpha$,
$\| (-\triangle)^\alpha f\|_{W^{s,p}(\Omega)}=\| f\|_{W^{s+2\alpha,p}(\Omega)}$ will not hold. By using Fourier transform, we have
$\| (-\triangle)^\alpha f\|_{W_\rho^{s,p}(\Omega)}=\| f\|_{W_\rho^{s+2\alpha,p}(\Omega)}$.

Finally, by \cite{applebaum2009}, let $X$ be a rotationally invariant stable process of index $\alpha\in (0,1)$. Its symbol is given
by $\eta(f)=-|f|^\alpha$ for all $f\in \mathbf{R}^n$. It is instructive to accept  that $\eta$ is the symbol for a legitimate
differential operator; then, using the usual correspondence $f_j\to -i\partial_j$ for $1\leq j\leq n$, we would write
\begin{eqnarray}A=\eta(\Omega)=-\left(\sqrt{-\partial_1^2-\partial_2^2-\cdots-\partial_n^2}\right)^\alpha=-(-\triangle )^{\alpha/2}.\end{eqnarray}
In fact, it is very useful to interpret $\eta(\Omega)$ as a fractional power of the operator $-\triangle$.
However, for the fractional Laplacian operator, it  can be defined as the generator of $\alpha$-stable L$\acute{e}$vy processes. More precisely,
if $X_t$ is the isotropic  $\alpha$-stable L$\acute{e}$vy processes starting at zero and $f$ is a smooth function, then
\begin{eqnarray}
(-\triangle)^{\alpha/2} f(x)= \lim\limits_{h\to 0^+} \frac{\mathbf{E}[f(x)-f(x+X_h)]}{h}.
\end{eqnarray}
This also indicates that the fractional power of operator $-\triangle$ and the fractional Laplacian operator are different.

\subsection{Fractional Laplacian operator and fractional derivative}

This subsection is focused on the relationship between the fractional Laplacian operator and the Riesz fractional derivative.

\begin{definition}\cite{gorenflo1998}
The Riesz fractional operator for $n-1<\alpha\leq n$ on a finite interval $0\leq x\leq l$ is defined as
\begin{eqnarray}
\frac{\partial^\alpha}{\partial |x|^\alpha}f(x) =\frac{-1}{2\cos(\frac{\alpha \pi}{2})}\left[{}_0D_x^\alpha+{}_xD_l^\alpha\right]f(x),
\end{eqnarray}
where
\begin{eqnarray*}
_0D_x^\alpha f(x)=\frac{1}{\Gamma(n-\alpha)}\frac{\partial^n}{\partial x^n}\int_0^x{(x-\eta)^{-\alpha-1+n}f(\eta)}d\eta
\end{eqnarray*}
and
\begin{eqnarray*}
_xD_l^\alpha f(x)=\frac{(-1)^n}{\Gamma(n-\alpha)}\frac{\partial^n}{\partial x^n}\int_x^l{(\eta-x)^{-\alpha-1+n}f(\eta)}d\eta
\end{eqnarray*}
are the left-sided and right-sided Riemann-Liouville fractional derivative, respectively.
\end{definition}

Moreover, according to \cite{samko1993}, the fractional Laplacian is the operator with symbol $|x|^\alpha.$  In other words, the following formula holds:
\begin{eqnarray}\label{fourierd}
(-\triangle )^{\alpha/2}f(x)= \mathcal{F}^{-1}|x|^\alpha \mathcal{F} f(x),
\end{eqnarray}
where $\mathcal{F}$ and $\mathcal{F}^{-1}$ denote the Fourier transform and inverse Fourier transform of $f(x)$, respectively.
We refer the readers  to \cite{landkof1972} for a detailed proof of the equivalence between the two definitions $(\ref{cauchyd})$
and $(\ref{fourierd})$ of fractional Laplacian operator.

By using Luchko's theorem in \cite{luchko1999},  we obtain the  following result on  the  equivalent
relationship between the Riesz fractional derivative $\frac{\partial^\alpha}{\partial |x|^\alpha}$ and
the fractional Laplacian operator $-(-\triangle)^{\alpha/2}$:
\begin{lemma}\cite{jliufw2012}\label{lemma2.1}
For a function $f(x)$ defined on the finite domain $[0,l]$, and $f(0)=f(l)=0$,  the following equality holds:
\begin{eqnarray}
-\left(-\triangle\right)^{\alpha/2} f(x)= \frac{-1}{2\cos(\frac{\alpha \pi}{2})}\left[{}_0D_x^\alpha f(x)+{}_xD_l^\alpha f(x)\right]=\frac{\partial^\alpha}{\partial |x|^\alpha}f(x),
\end{eqnarray}
where $\alpha\in (1,2)$ and the space fractional derivative $\frac{\partial^\alpha}{\partial |x|^\alpha}$ is a Riesz fractional derivative.
\end{lemma}

For more information on the analytical solution of the generalized multi-term time and space fractional partial
differential equations with Dirichlet nonhomogeneous boundary conditions, we refer the readers to \cite{jliufw2012}.
 For more information on the numerical solution
of fractional partial differential equation with Riesz space fractional derivatives on a finite domain, consult  \cite{liufw2010,liufw2004}.

\subsection{Fractional derivative and fractional power of operator}
In this part, we first show the following definition of the positive operator.
\begin{definition}\cite{ashyralyev2009}
The operator $A$ is said to be positive if its spectrum $\sigma(A)$ lies in the interior of the sector of angle $\varphi\in (0,\pi)$,
 symmetric with respect to the real axis, and if on the edges of this sector, $S_1=\{\rho e^{i\varphi}: 0\leq \rho <\infty\}$
and $S_2=\{\rho e^{-i\varphi}: 0\leq \rho <\infty\}$, and outside it the resolvent $(\lambda I-A)^{-1}$ is subject to the bound
\begin{eqnarray}
\left\|(\lambda I-A)^{-1}\right\|\leq \frac{M(\varphi)}{1+|\lambda|}.
\end{eqnarray}
\end{definition}
Moreover, for a positive operator $A$, any $\alpha>0$, one can define the negative fractional power of operator $A$ by the following formula \cite{komatsu1966-1}:
\begin{eqnarray}\label{5.2}
A^{-\alpha}=\frac{1}{2\pi i}\int_\Gamma{\lambda^{-\alpha} R(\lambda)}d\lambda,~~\left(R(\lambda)=(A-\lambda I)^{-1},~\Gamma=S_1 \cup S_2\right).
\end{eqnarray}
It is then quite easy to see that $A^{-\alpha}$ is a bounded operator, which is an entire function of $\alpha$, satisfying
$A^{-\alpha}=A^{-n}$ if $\alpha$ is an integer $n$, and $A^{-(\alpha+\beta)}=A^{-\alpha}A^{-\beta}$ for all $\alpha,\beta\in \mathbf{C}$.
Using $(\ref{5.2})$,  we have
\begin{eqnarray}
A^{-\alpha}=\frac{1}{2\pi i}\int_{-\infty}^0{\lambda^{-\alpha} R(\lambda)}d\lambda+\frac{1}{2\pi i}\int_0^{-\infty}{\lambda^{-\alpha} R(\lambda)}d\lambda,
\end{eqnarray}
Then taking the integration along the lower and upper sides of the cut respectively: $\lambda=se^{-\pi i}$ and $\lambda=se^{\pi i}$,  it follows that
\begin{eqnarray*}
A^{-\alpha}&=&\frac{e^{\alpha \pi i}}{2\pi i}\int_0^\infty{s^{-\alpha} R(-s)}ds- \frac{e^{-\alpha \pi i}}{2\pi i}\int_0^\infty{s^{-\alpha} R(-s)}ds\\
&=& \frac{\cos(\alpha \pi)+i\sin(\alpha \pi)}{2\pi i}\int_0^\infty{s^{-\alpha} R(-s)}ds- \frac{\cos(\alpha \pi)-i\sin(\alpha \pi)}{2\pi i}\int_0^\infty{s^{-\alpha} R(-s)}ds\\
&=& \frac{\sin(\alpha \pi)}{\pi}\int_0^\infty{s^{-\alpha} R(-s)}ds\\
&=&  \frac{1}{\Gamma(\alpha)\Gamma(1-\alpha)}\int_0^\infty{s^{-\alpha} R(-s)}ds
\end{eqnarray*}
Moreover, for any $\alpha\in (n-1,n)$, we get that
\begin{eqnarray}\label{2.18}
A^{\alpha}f=A^{\alpha-n}A^{n}f=  \frac{1}{\Gamma(n-\alpha)\Gamma(1+\alpha-n)}\int_0^\infty{s^{\alpha-n} R(-s)A^{n}f}ds.
\end{eqnarray}
and for more properties on the fractional power of a positive operator, please see
\cite{komatsu1967-II,komatsu1969-III,carracedo2001,balakrishnan1960} and the references cited therein.

Now we are ready to state the following results on the connection of fractional derivative and integral with fractional power of positive operator.

\begin{theorem}
Let the absolutely space \[AC^n[0,l]: =\{f: f^{(n-1)}(x) \in C[0, l],~ f^{(n)}(x) \in L^2[0, l]\}\] and
let $A$ be the operator defined by the formula $Af(x) = f'(x)$ with the domain $\{f: f\in AC^n[0,l], ~f^{(n)}(0)=0\}$.
Then $A$ is a positive operator in the Banach space $AC^n[0,l]$ and
\begin{equation}A^\alpha f(x)={}_0D_x^\alpha f(x),~~n-1<\alpha<n\end{equation}
for all $f (x) \in D(A)$.
\end{theorem}
\textbf{Proof.} By \cite{FODO}, the operator $A+sI$ $(s\geq 0)$ has a bounded inverse and the resolvent of $A$ is given by
\begin{equation}
\left(\left(A+sI\right)^{-1}f\right)(x)=\int_0^x{e^{-s(x-y)}f(y)}dy.
\end{equation}
Then the operator $A$ is a positive operator in $AC^n[0,l]$ and  Eq. $(\ref{2.18})$ gives
\begin{eqnarray*}
A^{\alpha}f(x)&=& \frac{1}{\Gamma(n-\alpha)\Gamma(1+\alpha-n)}\int_0^\infty{s^{\alpha-n} R(-s)A^{n}f(x)}ds\\
&=&\frac{1}{\Gamma(n-\alpha)\Gamma(1+\alpha-n)}\int_0^\infty{s^{\alpha-n} \left(A+sI\right)^{-1}f^{(n)}(x)}ds\\
&=& \frac{1}{\Gamma(n-\alpha)\Gamma(1+\alpha-n)}\int_0^\infty{s^{\alpha-n} \int_0^x{e^{-s(x-y)}f^{(n)}(y)}dy}ds\\
&=& \frac{1}{\Gamma(n-\alpha)\Gamma(1+\alpha-n)} \int_0^x\left[\int_0^\infty{s^{\alpha-n}e^{-s(x-y)}}ds\right]f^{(n)}(y)dy.
\end{eqnarray*}
Let $s(x-y)=\lambda$. we get that
\begin{eqnarray*}
\int_0^\infty{s^{\alpha-n}e^{-s(x-y)}}ds=(x-y)^{n-\alpha-1}\int_0^\infty{\lambda^{\alpha-n}e^{-\lambda}}d\lambda=(x-y)^{n-\alpha-1}\Gamma(\alpha-n+1),
\end{eqnarray*}
Then
\begin{eqnarray*}
A^{\alpha}f(x)&=& \frac{1}{\Gamma(n-\alpha)\Gamma(1+\alpha-n)} \int_0^x{(x-y)^{n-\alpha-1}\Gamma(\alpha-n+1)f^{(n)}(y)}dy\\
&=&\frac{1}{\Gamma(n-\alpha)} \int_0^x{(x-y)^{n-\alpha-1}f^{(n)}(y)}dy  \\
&=& {}_0D_x^{\alpha} f(x).
\end{eqnarray*}
This completes the proof.

\section{The emerging research opportunities}
\label{sec3}

\setcounter{section}{3}\setcounter{equation}{0}
Recent advances in modeling and control of fractional diffusion systems, fractional reaction-diffusion systems
and fractional reaction-diffusion-advection systems have been reviewed in the framework of CPSs. The fractional order
DPSs have now been found wide applications for describing many physical phenomena, such as sub-diffusion or super-diffusion processes.
At the same time, to our best knowledge, many problems are still open calling for research cooperation of multi-disciplines such as mathematical
modelling, engineering applications, and information sciences.

First, it is worth noting that in the more recent monograph \cite{umarov2015}, the theory of pseudo-differential operators with
singular symbols,  and the connections between them and those three types of operators are explored. See
\cite{umarov2012,umarov1991,umarov2011} for more knowledge on pseudo-differential operator. Moreover, we claim that those
equivalences between fractional Laplacian operator and fractional derivative, fractional order of operator and fractional derivative
discussed in this paper can introduce new mathematical vehicles to study   fractional order generalized DPSs. For example, when we
study a fractional DPSs with Riesz fractional derivative, by Lemma $\ref{lemma2.1}$, the spectral representation methods can be
used to characterize the solution of the dynamic system. Then we can study the existence of solutions,
stability, controllability and observability of the system under consideration.

Potential topics such as modeling the sub-diffusion or super-diffusion processes with consideration of the networked mobile actuators and mobile sensors, the communication among the actuators and sensors, collocated  or non-collocated actuators and sensors, their robustness and optimality
problems are all interesting and worthy  much more efforts in future. Another interesting and important topic   concerns  the time-space
fractional DPSs  where the traditional first order derivative is replaced by a fractional order derivative with  respect to the time $t$.

Furthermore, in the case of diffusion systems, it is worth mentioning that, in general, not all the states can be reached in the whole domain of interest
\cite{AEIJAI,afifi2002regionally,ge2015arxiv} and it would be more challenging in nature since the dynamics of the real-life control problem is always
hybrid continuous and discrete. Due to the spatial-temporal sampling and discrete nature of decision and control, the notions of regional analysis should be
introduced, i.e., we can consider the regional stability, regional controllability, regional observability etc. of the system under consideration.
In addition, as stated in \cite{chen2004}, from an application point of view, some plain questions such as ``How many actuators/sensors are sufficient
and how to best configure them for a fractional DPSs control process?", ``Given the desirable zone shape, is it possible to control or contain the fractional diffusion
process within the given zone?", if not, ``how to quantify the controllability/observability of the actuators/sensors" and
etc. might be asked, which in fact raises some important theoretical challenges and open new opportunities for further research.

\section{Conclusion}
This paper is concerned with the fractional order DPSs with three different operators:
fractional Laplacian operator, fractional power of operator  and fractional derivative.
The relationship among the three operators and the emerging research opportunities are
introduced. We hope that the results here could provided some insight into the control
theory analysis of fractional order DPSs in particular and CPSs in general. The results
presented here can also be extended to complex fractional order DPSs and  various open
questions are still pending.  For instance, the problem of  regional optimal control of
fractional order DPSs with more complicated sensing and actuation configurations are of
great interest.

\noindent {\bf Acknowledgement.}
This work was supported  by Chinese Universities Scientific Fund (No.CUSF-DH-D-2014061)
 and  Natural Science Foundation of Shanghai (No.15ZR1400800).

\end{document}